\documentclass[11pt]{amsart}
\usepackage{amsmath,amssymb,mathrsfs}
\newtheorem{theorem}{Theorem}
\newtheorem{proposition}[theorem]{Proposition}
\newtheorem{corollary}[theorem]{Corollary}
\newtheorem{lemma}[theorem]{Lemma}
\newcommand{\tracefreeRic}{\overset{\text{\rm o}}{\text{\rm Ric}}}

\begin{document}

\title{Area-minimizing projective planes in three-manifolds}
\author{H. Bray, S. Brendle, M. Eichmair, and A. Neves}
\address{Department of Mathematics \\ Duke University \\ Durham, NC 27708}
\address{Department of Mathematics \\ Stanford University \\ 450 Serra Mall, Bldg 380 \\ Stanford, CA 94305} 
\address{Department of Mathematics \\ MIT \\ 77 Massachusetts Avenue \\ Cambridge, MA 02139}
\address{Imperial College \\ Huxley Building \\ 180 Queen's Gate \\ London SW7 2RH \\ United Kingdom}



\maketitle

\section{Introduction}

Let $M$ be a compact three-manifold equipped with a Riemannian metric $g$. We denote by $\mathscr{F}$ the set of all embedded surfaces $\Sigma \subset M$ such that $\Sigma$ is homeomorphic to $\mathbb{RP}^2$.

Throughout this paper, we shall assume that $\mathscr{F}$ is non-empty. We define 
\begin{equation} 
\label{definition.of.A}
\mathscr{A}(M,g) = \inf \{\text{\rm area}(\Sigma,g): \Sigma \in \mathscr{F}\}. 
\end{equation}
Recall that the systole of $(M,g)$ is defined as 
\begin{equation} 
\label{definition.of.systole}
\text{\rm sys}(M,g) = \inf \{L(\gamma): \text{$\gamma$ is a non-contractible loop in $M$}\} 
\end{equation}
(see e.g. \cite{Gromov}). The definition (\ref{definition.of.A}) is similar in spirit to (\ref{definition.of.systole}). Rather than minimizing lengths of non-contractible loops, we minimize area among embedded projective planes. The quantity $\mathscr{A}(M,g)$ is also related to the notion of width studied by Colding and Minicozzi \cite{Colding-Minicozzi1}, \cite{Colding-Minicozzi2}.

Our main result is the following:

\begin{theorem}
\label{theorem.1}
Let $(M,g)$ be a compact three-manifold equipped with a Riemannian metric. Moreover, we assume that $M$ contains an embedded projective plane. Then 
\begin{equation}
\label{upper.bound.for.A}
\mathscr{A}(M,g) \, \inf_M R_g \leq 12\pi 
\end{equation} 
and 
\begin{equation} 
\label{lower.bound.for.A}
\mathscr{A}(M,g) \geq \frac{2}{\pi} \, \text{\rm sys}(M,g)^2 > 0. 
\end{equation} 
Here, $R_g$ denotes the scalar curvature of the metric $g$.
\end{theorem}

Combining (\ref{upper.bound.for.A}) and (\ref{lower.bound.for.A}) yields 
\[\text{\rm sys}(M,g)^2 \, \inf_M R_g \leq 6\pi^2.\] 
We note that Gromov and Lawson proved that 
\[\text{\rm Rad}(M,g)^2 \, \inf_M R_g \leq 4\pi^2,\] 
where $\text{\rm Rad}(M,g)$ denotes the homology filling radius of $(M,g)$ (see \cite{Gromov-Lawson}, Theorem $\text{\rm G}_2$). A similar result was established by Schoen and Yau (cf. \cite{Schoen-Yau}, Theorem 1).

The inequalities (\ref{upper.bound.for.A}) and (\ref{lower.bound.for.A}) are both sharp on $\mathbb{RP}^3$. Indeed, if $g$ denotes the round metric on $\mathbb{RP}^3$ with constant sectional curvature $1$, then $R_g = 6$ and $\text{\rm sys}(\mathbb{RP}^3,g) = \pi$. Using (\ref{upper.bound.for.A}) and (\ref{lower.bound.for.A}), we conclude that $\mathscr{A}(\mathbb{RP}^3,g) = 2\pi$.

We next characterize the case of equality in (\ref{upper.bound.for.A}).

\begin{theorem}
\label{theorem.2}
Let $(M,g)$ be a compact three-manifold equipped with a Riemannian metric. Moreover, we assume that $M$ contains an embedded projective plane. If $\mathscr{A}(M,g) \, \inf_M R_g = 12\pi$, then $(M,g)$ is isometric to $\mathbb{RP}^3$ up to scaling.
\end{theorem}

In particular, if $\text{\rm sys}(M,g)^2 \, \inf_M R_g = 6\pi^2$, then $(M,g)$ is isometric to $\mathbb{RP}^3$ up to scaling.

We now describe the proof of Theorem \ref{theorem.1}. The inequality (\ref{lower.bound.for.A}) follows directly from a classical theorem due to Pu \cite{Pu}. The proof of (\ref{upper.bound.for.A}) is more subtle. General results of Meeks, Simon, and Yau \cite{Meeks-Simon-Yau} imply that the infimum in (\ref{definition.of.A}) is attained by an embedded surface $\Sigma \in \mathscr{F}$. The inequality (\ref{lower.bound.for.A}) is then obtained using special choices of variations in the second variation formula. When $\Sigma$ is two-sided, we consider unit-speed variations. When $\Sigma$ is one-sided, we use a technique due to Hersch \cite{Hersch} to construct suitable sections of the normal bundle. This trick has also been used in other contexts, see e.g. \cite{Colding-Minicozzi1}, \cite{Huisken-Yau}, \cite{Li-Yau}.

In order to prove Theorem \ref{theorem.2}, we assume that $g_0$ is a Riemannian metric on $M$ satisfying $\mathscr{A}(M,g_0) \, \inf_M R_{g_0} = 12\pi$. By scaling, we may assume that $\mathscr{A}(M,g_0) = 2\pi$ and $\inf_M R_{g_0} = 6$. We then evolve the metric $g_0$ by Hamilton's Ricci flow. We show that $\mathscr{A}(M,g(t)) \geq 2\pi (1-4t)$ and $\inf_M R_{g(t)} \geq \frac{6}{1-4t}$. Using Theorem \ref{theorem.1}, we conclude that both inequalities are, in fact, equalities. The strict maximum principle then implies that $(M,g(t))$ has constant sectional curvature for each $t$.

We thank Fernando Marques for discussions, and Ian Agol and Christina Sormani for comments on an earlier version of this paper.

\section{Proof of Theorem \ref{theorem.1}}

In this section, we present the proof of Theorem \ref{theorem.1}. As above, we assume that $M$ is a compact three-manifold which contains an embedded projective plane. In order to verify (\ref{lower.bound.for.A}), we need the following result: 

\begin{proposition}
\label{incompressibility}
Let $\Sigma$ be an arbitrary surface in $\mathscr{F}$. Then the induced map $i_\#: \pi_1(\Sigma) \to \pi_1(M)$ is injective.
\end{proposition}

\textbf{Proof.} 
We argue by contradiction. If the map $i_\#: \pi_1(\Sigma) \to \pi_1(M)$ fails to be injective, then the bundle $TM|_\Sigma$ is orientable. Since the tangent bundle $T\Sigma$ is non-orientable, we conclude that $\Sigma$ has non-trivial normal bundle. Let $\gamma: [0,1] \to \Sigma$ be a smooth closed curve in $\Sigma$ which represents a non-trivial element of $\pi_1(\Sigma)$. For each $t \in [0,1]$, we can find a unit vector $\nu(t) \in T_{\gamma(t)} M$ which is orthogonal to the tangent space $T_{\gamma(t)} \Sigma$. Moreover, we may assume that $\nu(t)$ depends continuously on $t$. Since $\Sigma$ has non-trivial normal bundle, we have $\nu(0) = -\nu(1)$. 

For each $\varepsilon > 0$, we define a path $\gamma_\varepsilon: [0,1] \to M$ by 
\[\gamma_\varepsilon(t) = \exp_{\gamma(t)}(\varepsilon \, \sin(\pi t) \, \nu(t)).\] 
Clearly, $\gamma_\varepsilon$ is a smooth closed curve in $M$. If we choose $\varepsilon > 0$ sufficiently small, then the curve $\gamma_\varepsilon$ intersects $\Sigma$ in exactly one point, and the intersection is transversal. Consequently, $\gamma_\varepsilon$ represents a non-trivial element of $\pi_1(M)$. Since $\gamma_\varepsilon$ is homotopic to $\gamma$, it follows that $\gamma$ represents a non-trivial element of $\pi_1(M)$. This is a contradiction. \\

Combining Proposition \ref{incompressibility} with Pu's inequality, we can draw the following conclusion:

\begin{corollary}
\label{A.sys}
We have $\mathscr{A}(M,g) \geq \frac{2}{\pi} \, \text{\rm sys}(M,g)^2 > 0$. 
\end{corollary}

\textbf{Proof.} 
Fix an arbitrary surface $\Sigma \in \mathscr{F}$. Then $\Sigma$ is homeomorphic to $\mathbb{RP}^2$, and the induced map $i_\#: \pi_1(\Sigma) \to \pi_1(M)$ is injective. Using Pu's inequality (Theorem 1 in \cite{Pu}), we obtain 
\[\text{\rm area}(\Sigma,g) \geq \frac{2}{\pi} \, \text{\rm sys}(\Sigma,g)^2 \geq \frac{2}{\pi} \, \text{\rm sys}(M,g)^2.\] 
Since $\Sigma \in \mathscr{F}$ is arbitrary, the assertion follows. \\

We now describe the proof of (\ref{upper.bound.for.A}). In the first step, we show that the infimum in (\ref{definition.of.A}) is attained by some surface $\Sigma \in \mathscr{F}$. To that end, we employ a general theorem of Meeks, Simon, and Yau \cite{Meeks-Simon-Yau} (see also \cite{Hass-Scott}, Theorem 5.2).

\begin{proposition}
\label{existence.of.minimizer}
There exists a surface $\Sigma \in \mathscr{F}$ such that $\text{\rm area}(\Sigma,g) = \mathscr{A}(M,g)$.
\end{proposition}

\textbf{Proof.} 
We can find a sequence of surfaces $\Sigma_k \in \mathscr{F}$ such that 
\[\text{\rm area}(\Sigma_k,g) \leq \mathscr{A}(M,g) + \varepsilon_k,\] 
where $\varepsilon_k \to 0$ as $k \to \infty$. This implies 
\[\text{\rm area}(\Sigma_k,g) \leq \inf_{\Sigma \in \mathscr{J}(\Sigma_k)} \text{\rm area}(\Sigma,g) + \varepsilon_k,\] 
where $\mathscr{J}(\Sigma_k)$ denotes the collection of all embedded surfaces isotopic to $\Sigma_k$. By 
Theorem 1 in \cite{Meeks-Simon-Yau}, a subsequence of the sequence $\Sigma_k$ converges weakly to a disjoint union of smooth embedded minimal surfaces $\Sigma^{(1)}, \hdots, \Sigma^{(R)}$ with positive integer multiplicities. More precisely, we can find positive integers $R,n_1,\hdots,n_R$ and pairwise disjoint embedded minimal surfaces $\Sigma^{(1)}, \hdots, \Sigma^{(R)}$ such that 
\[\sum_{j=1}^R n_j \int_{\Sigma^{(j)}} f \, d\mu_g = \lim_{k \to \infty} \int_{\Sigma_k} f \, d\mu_g\] 
for every continuous function $f: M \to \mathbb{R}$. In particular, we have 
\begin{equation} 
\label{area}
\sum_{j=1}^R n_j \, \text{\rm area}(\Sigma^{(j)},g) \leq \mathscr{A}(M,g). 
\end{equation}
Following Meeks, Simon, and Yau \cite{Meeks-Simon-Yau}, we define surfaces $S_k^{(1)}, \hdots, S_k^{(R)}$ as follows: if $n_j = 2m_j$ is even, then $S_k^{(j)}$ is defined by 
\[S_k^{(j)} = \bigcup_{r=1}^{m_j} \Big \{ x \in M: d(x,\Sigma^{(j)}) = \frac{r}{k} \Big \}\] 
On the other hand, if $n_j = 2m_j + 1$ is odd, then $S_k^{(j)}$ is defined by 
\[S_k^{(j)} = \Sigma^{(j)} \cup \bigcup_{r=1}^{m_j} \Big \{ x \in M: d(x,\Sigma^{(j)}) = \frac{r}{k} \Big \}.\] 
By Remark 3.27 in \cite{Meeks-Simon-Yau}, we can find embedded surfaces $S_k^{(0)}$ and $\tilde{\Sigma}_k$ with the following properties: 
\begin{itemize}
\item[(i)] The surface $S_k = \bigcup_{j=0}^R S_k^{(j)}$ is isotopic to $\tilde{\Sigma}_k$ if $k$ is sufficiently large.
\item[(ii)] The surface $\tilde{\Sigma}_k$ is obtained from $\Sigma_{q_k}$ by $\gamma_0$-reduction (cf. \cite{Meeks-Simon-Yau}, Section 3). 
\item[(iii)] We have $S_k^{(0)} \cap \big ( \bigcup_{j=1}^R S_k^{(j)} \big ) = \emptyset$. Moreover, $\text{\rm area}(S_k^{(0)},g) \to 0$ as $k \to \infty$.
\end{itemize}
By assumption, $\Sigma_{q_k}$ is homeomorphic to $\mathbb{RP}^2$, and $\tilde{\Sigma}_k$ is obtained from $\Sigma_{q_k}$ by $\gamma_0$-reduction. Consequently, one of the connected components of $\tilde{\Sigma}_k$ is homeomorphic to $\mathbb{RP}^2$. 

Hence, if $k$ is sufficiently large, then one of the connected components of $S_k$ is homeomorphic to $\mathbb{RP}^2$. Let us denote this connected component by $E_k$. Since $E_k \in \mathscr{F}$, we have $\text{\rm area}(E_k,g) \geq \mathscr{A}(M,g) > 0$. On the other hand, we have $\text{\rm area}(S_k^{(0)},g) \to 0$ as $k \to \infty$. This implies $\text{\rm area}(E_k,g) > \text{\rm area}(S_k^{(0)},g)$ if $k$ is sufficiently large. Hence, if $k$ is sufficiently large, then $E_k$ cannot be contained in $S_k^{(0)}$. Since $E_k \subset S_k$ is connected, it follows that $E_k$ is a connected component of $S_k^{(i)}$ for some integer $i \in \{1,\hdots,R\}$. Thus, $E_k$ is either homeomorphic to $\Sigma^{(i)}$ or to a double cover of $\Sigma^{(i)}$. Since $E_k$ is homeomorphic to $\mathbb{RP}^2$, we conclude that $\Sigma^{(i)}$ is homeomorphic to $\mathbb{RP}^2$. This shows that $\Sigma^{(i)} \in \mathscr{F}$. Moreover, it follows from (\ref{area}) that $\text{\rm area}(\Sigma^{(i)},g) \leq \mathscr{A}(M,g)$. Therefore, the surface $\Sigma^{(i)}$ is the desired minimizer. \\

\begin{proposition} 
\label{stability} 
Let $\Sigma$ be a surface in $\mathscr{F}$ satisfying $\text{\rm area}(\Sigma,g) = \mathscr{A}(M,g)$. Then 
\[\int_\Sigma (\text{\rm Ric}_g(\nu,\nu) + |I\!I|^2) \, d\mu_g \leq 4\pi.\] 
\end{proposition}

\textbf{Proof.} 
By the uniformization theorem, we can find a diffeomorphism $\varphi: \mathbb{RP}^2 \to \Sigma$ such that the metric $\varphi^* g$ is conformal to the standard metric on $\mathbb{RP}^2$. We may lift the map $\varphi: \mathbb{RP}^2 \to \Sigma$ to a map $\hat{\varphi}: S^2 \to \Sigma$. Clearly, $\hat{\varphi}(x) = \hat{\varphi}(-x)$ for all $x \in S^2$. Moreover, the metric $\hat{\varphi}^* g$ is conformal to the standard metric $h$ on $S^2$.

We next consider the pull-back of the normal bundle $N\Sigma$ under the map $\hat{\varphi}: S^2 \to \Sigma$. Since the bundle $\hat{\varphi}^* N\Sigma$ is trivial, we can find a smooth section $\nu \in \Gamma(\hat{\varphi}^* N\Sigma)$ such that $|\nu(x)| = 1$ for all $x \in S^2$. For each point $x \in S^2$, the vector $\nu(x)$ is a unit normal vector to $\Sigma$ at the point $\hat{\varphi}(x)$. There are two possibilities:

\textit{Case 1:} Suppose that $\Sigma$ is two-sided, so that $\nu(x) = \nu(-x)$ for all $x \in S^2$. In this case, there exists a section $V \in \Gamma(N\Sigma)$ such that $\nu(x) = V(\hat{\varphi}(x))$ for all $x \in S^2$. Since $\Sigma$ has minimal area among all surfaces in $\mathscr{F}$, we have 
\[\int_\Sigma (\text{\rm Ric}_g(\nu,\nu) + |I\!I|^2) \, d\mu_g \leq \int_\Sigma |\nabla V|^2 \, d\mu_g = 0.\] 

\textit{Case 2:} We now assume that $\Sigma$ is one-sided, so that $\nu(x) = -\nu(-x)$ for all $x \in S^2$. We may identify $S^2$ with the unit sphere in $\mathbb{R}^3$. For each $j \in \{1,2,3\}$, we define a section $\sigma_j \in \Gamma(\hat{\varphi}^* N\Sigma)$ by $\sigma_j(x) = x_j \, \nu(x)$ for all $x \in S^2$. Note that $\sigma_j(x) = \sigma_j(-x)$ for all $x \in S^2$. Hence, there exists a section $V_j \in \Gamma(N\Sigma)$ such that $\sigma_j(x) = V_j(\hat{\varphi}(x))$ for all $x \in S^2$. Since $\sum_{j=1}^3 |\sigma_j(x)|^2 = 1$ for all $x \in S^2$, we conclude that $\sum_{j=1}^3 |V_j|^2 = 1$ at each point on $\Sigma$. 

Since $\Sigma$ has minimal area among all surfaces in $\mathscr{F}$, we have 
\[\int_\Sigma (\text{\rm Ric}_g(\nu,\nu) + |I\!I|^2) \, |V_j|^2 \, d\mu_g \leq \int_\Sigma |\nabla V_j|^2 \, d\mu_g\] 
for each $j \in \{1,2,3\}$. Since the metric $\hat{\varphi}^* g$ is conformal to the standard metric $h$ on $S^2$, we have 
\[\int_\Sigma |\nabla V_j|^2 \, d\mu_g = \frac{1}{2} \int_{S^2} |\nabla x_j|_{\hat{\varphi}^* g}^2 \, d\mu_{\hat{\varphi}^* g} = \frac{1}{2} \int_{S^2} |\nabla x_j|_h^2 \, d\mu_h\] 
for each $j \in \{1,2,3\}$. Using the identity $\Delta_h x_j + 2x_j = 0$, we conclude that 
\[\int_\Sigma (\text{\rm Ric}_g(\nu,\nu) + |I\!I|^2) \, |V_j|^2 \, d\mu_g \leq \frac{1}{2} \int_{S^2} |\nabla x_j|_h^2 \, d\mu_h = \int_{S^2} x_j^2 \, d\mu_h = \frac{4\pi}{3}\] 
for each $j \in \{1,2,3\}$. Summation over $j$ yields
\[\int_\Sigma (\text{\rm Ric}_g(\nu,\nu) + |I\!I|^2) \, d\mu_g \leq 4\pi,\] 
as claimed. \\

\begin{proposition}
\label{gauss.bonnet}
Let $\Sigma$ be an arbitrary surface in $\mathscr{F}$. Then 
\[\int_\Sigma (R_g - 2 \, \text{\rm Ric}_g(\nu,\nu) - |I\!I|^2) \, d\mu_g \leq 4\pi.\] 
\end{proposition}

\textbf{Proof.} 
Using the Gauss equation, we obtain 
\[R_g - 2 \, \text{\rm Ric}_g(\nu,\nu) - |I\!I|^2 = 2K - |H|^2,\] 
where $K$ is the Gaussian curvature of $\Sigma$ and $H$ denotes its mean curvature vector. Since $\Sigma$ is homeomorphic to $\mathbb{RP}^2$, we conclude that 
\[\int_\Sigma (R_g - 2 \, \text{\rm Ric}_g(\nu,\nu) - |I\!I|^2) \, d\mu_g \leq 2 \int_\Sigma K \, d\mu_g = 4\pi\] 
by the Gauss-Bonnet theorem. \\

\begin{corollary}
We have $\mathscr{A}(M,g) \, \inf_M R_g \leq 12\pi$. 
\end{corollary}

\textbf{Proof.} 
By Proposition \ref{existence.of.minimizer}, there exists an embedded surface $\Sigma \in \mathscr{F}$ such that $\text{\rm area}(\Sigma,g) = \mathscr{A}(M,g)$. Using Proposition \ref{stability} and Proposition \ref{gauss.bonnet}, we obtain 
\begin{align*} 
\mathscr{A}(M,g) \, \inf_M R_g 
&= \text{\rm area}(\Sigma,g) \, \inf_M R_g \\ 
&\leq \int_\Sigma (R_g + |I\!I|^2) \, d\mu_g \\ 
&\leq 4\pi + 2 \int_\Sigma (\text{\rm Ric}_g(\nu,\nu) + |I\!I|^2) \, d\mu_g \\ 
&\leq 12\pi. 
\end{align*}
This completes the proof. \\

\section{Proof of Theorem \ref{theorem.2}}

In this section, we analyze the case of equality in (\ref{upper.bound.for.A}). To that end, we fix a Riemannian metric $g_0$ on $M$. By a theorem of Hamilton \cite{Hamilton1}, there exists a real number $T > 0$ and a family of metrics $g(t)$, $t \in [0,T]$, such that $g(0) = g_0$ and 
\begin{equation} 
\label{ricci.flow}
\frac{\partial}{\partial t} g(t) = -2 \, \text{\rm Ric}_{g(t)} 
\end{equation} 
for all $t \in [0,T]$ (see also \cite{DeTurck}). The evolution equation (\ref{ricci.flow}) is known as the Ricci flow, and plays an important role in Riemannian geometry (see e.g. \cite{Brendle}, \cite{Brendle-Schoen}, \cite{Hamilton1}, \cite{Hamilton2}).

\begin{lemma}
The function $t \mapsto \mathscr{A}(M,g(t))$ is Lipschitz continuous.
\end{lemma}

\textbf{Proof.} 
We can find a real number $\Lambda > 0$ such that $\sup_M |\text{\rm Ric}_{g(t)}| \leq \Lambda$ for all $t \in [0,T]$. This implies
\[e^{-2\Lambda |t_0 - t_1|} \, g(t_0) \leq g(t_1) \leq e^{2\Lambda |t_0 - t_1|} \, g(t_0)\] 
for all times $t_0,t_1 \in [0,T]$. Consequently, we have 
\[e^{-2\Lambda |t_0 - t_1|} \, \mathscr{A}(M,g(t_0)) \leq \mathscr{A}(M,g(t_1)) \leq e^{2\Lambda |t_0 - t_1|} \, \mathscr{A}(M,g(t_0))\] 
for all times $t_0,t_1 \in [0,T]$. From this, the assertion follows. \\

In the next step, we show that the function $\mathscr{A}(M,g(t)) + 8\pi t$ is increasing in $t$. This result is similar in spirit to a theorem of Hamilton regarding the evolution of the area of stable minimal two-spheres under the Ricci flow (see \cite{Hamilton-survey}, Section 12). Related results for curves can be found in \cite{Hamilton-survey} and \cite{Ilmanen-Knopf}.

\begin{proposition} 
\label{evolution.of.A}
We have 
\[\mathscr{A}(M,g(t)) \geq \mathscr{A}(M,g_0) - 8\pi t\] 
for all $t \in [0,T]$. 
\end{proposition}

\textbf{Proof.} 
Suppose the assertion is false. Then there exists a time $\tau \in (0,T]$ such that 
\[\mathscr{A}(M,g(\tau)) < \mathscr{A}(M,g_0) - 8\pi \tau.\] 
Hence, we can find a real number $\varepsilon > 0$ such that 
\[\mathscr{A}(M,g(\tau)) < \mathscr{A}(M,g_0) - 8\pi \tau - 2\varepsilon \tau.\] 
We next define 
\[t_0 = \inf \big \{ t \in [0,T]: \mathscr{A}(M,g(t)) < \mathscr{A}(M,g_0) - (8\pi + \varepsilon) t - \varepsilon \tau \big \}.\] 
Clearly, $t_0 \in (0,\tau)$. Moreover, we have 
\[\mathscr{A}(M,g(t_0)) - \mathscr{A}(M,g(t)) \leq -(8\pi + \varepsilon) \, (t_0 - t)\] 
for all $t \in [0,t_0)$. By Proposition \ref{existence.of.minimizer}, we can find an embedded surface $\Sigma \in \mathscr{F}$ satisfying 
\[\text{\rm area}(\Sigma,g(t_0)) = \mathscr{A}(M,g(t_0)).\] 
For this choice of $\Sigma$, we have 
\begin{align*} 
\text{\rm area}(\Sigma,g(t_0)) - \text{\rm area}(\Sigma,g(t)) 
&\leq \mathscr{A}(M,g(t_0)) - \mathscr{A}(M,g(t)) \\ 
&\leq -(8\pi + \varepsilon) \, (t_0 - t) 
\end{align*} 
for all $t \in [0,t_0)$. This implies 
\[\frac{d}{dt} \text{\rm area}(\Sigma,g(t)) \Big |_{t=t_0} \leq -8\pi - \varepsilon.\] 
On the other hand, it follows from (\ref{ricci.flow}) that 
\[\frac{d}{dt} \text{\rm area}(\Sigma,g(t)) \Big |_{t=t_0} = -\int_\Sigma (\text{\rm Ric}_{g(t_0)}(e_1,e_1) + \text{\rm Ric}_{g(t_0)}(e_2,e_2)) \, d\mu_{g(t_0)},\] 
where $\{e_1,e_2\}$ denotes a local orthonormal frame on $\Sigma$ with respect to the metric $g(t_0)$. Using Proposition \ref{stability} and Proposition \ref{gauss.bonnet}, we obtain  
\begin{align*} 
\frac{d}{dt} \text{\rm area}(\Sigma,g(t)) \Big |_{t=t_0} 
&= -\int_\Sigma (R_{g(t_0)} - \text{\rm Ric}_{g(t_0)}(\nu,\nu)) \, d\mu_{g(t_0)} \\ 
&\geq -4\pi - \int_\Sigma (\text{\rm Ric}_{g(t_0)}(\nu,\nu) + |I\!I|^2) \, d\mu_{g(t_0)} \\ 
&\geq -8\pi. 
\end{align*} 
This is a contradiction. \\

\begin{proposition}
\label{upper.bound.for.scal} 
Suppose that $\mathscr{A}(M,g_0) = 2\pi$. Then 
\[\inf_M R_{g(t)} \leq \frac{6}{1-4t}\] 
for all $t \in [0,T] \cap [0,\frac{1}{4})$. 
\end{proposition} 

\textbf{Proof.} 
By Theorem \ref{theorem.1}, we have 
\[\mathscr{A}(M,g(t)) \, \inf_M R_{g(t)} \leq 12\pi\] 
for all $t \in [0,T]$. Moreover, it follows from Proposition \ref{evolution.of.A} that 
\[\mathscr{A}(M,g(t)) \geq \mathscr{A}(M,g(0)) - 8\pi t = 2\pi \, (1-4t)\] 
for all $t \in [0,T]$. Putting these facts together, the assertion follows. \\

\begin{proposition}
\label{round}
Suppose that $\mathscr{A}(M,g_0) \, \inf_M R_{g_0} = 12\pi$. Then the manifold $(M,g_0)$ has constant sectional curvature.
\end{proposition}

\textbf{Proof.} 
After rescaling the metric if necessary, we may assume that $\mathscr{A}(M,g_0) = 2\pi$ and $\inf_M R_{g_0} = 6$. The scalar curvature of $g(t)$ satisfies the evolution equation 
\[\frac{\partial}{\partial t} R_{g(t)} = \Delta R_{g(t)} + 2 \, |\text{\rm Ric}_{g(t)}|^2.\] 
This identity can be rewritten as  
\[\frac{\partial}{\partial t} R_{g(t)} = \Delta R_{g(t)} + \frac{2}{3} \, R_{g(t)}^2 + 2 \, |\tracefreeRic_{g(t)}|^2,\] 
where $\tracefreeRic_{g(t)}$ denotes the trace-free Ricci tensor of $g(t)$. Using the maximum principle, we conclude that $T < \frac{1}{4}$ and 
\begin{equation} 
\label{lower.bound.for.scal}
\inf_M R_{g(t)} \geq \frac{6}{1-4t} 
\end{equation}
for all $t \in [0,T]$ (see e.g. \cite{Brendle-book}, Proposition 2.19). By Proposition \ref{upper.bound.for.scal}, the inequality (\ref{lower.bound.for.scal}) is an equality. Using the strict maximum principle, we obtain 
\[R_{g(t)} = \frac{6}{1-4t}\] 
on $M \times [0,T]$. Substituting this into the evolution equation for the scalar curvature, we deduce that $|\tracefreeRic_{g(t)}|^2 = 0$ on $M \times [0,T]$. Since the Weyl tensor vanishes in dimension $3$, it follows that $(M,g(t))$ has constant sectional curvature for all $t \in [0,T]$. \\

By Proposition \ref{round}, the universal cover of $(M,g_0)$ is isometric to $S^3$ up to scaling. Hence, it remains to analyze the fundamental group of $M$.

\begin{proposition}
Suppose that $\mathscr{A}(M,g_0) \, \inf_M R_{g_0} = 12\pi$. Then $|\pi_1(M)| = 2$.
\end{proposition}

\textbf{Proof.} 
By scaling, we may assume that $\mathscr{A}(M,g_0) = 2\pi$ and $\inf_M R_{g_0} = 6$. By Proposition \ref{existence.of.minimizer}, there exists a surface $\Sigma \in \mathscr{F}$ such that $\text{\rm area}(\Sigma,g_0) = \mathscr{A}(M,g_0)$. Using Proposition \ref{stability} and Proposition \ref{gauss.bonnet}, we obtain 
\[12\pi = \text{\rm area}(\Sigma,g_0) \, \inf_M R_{g_0} \leq \int_\Sigma (R_{g_0} + |I\!I|^2) \, d\mu_{g_0} \leq 12\pi.\] 
Therefore, the surface $\Sigma$ is totally geodesic. 

By Proposition \ref{round}, there exists a local isometry $F: S^3 \to (M,g_0)$. Note that $F$ is a covering map (cf. \cite{Cheeger-Ebin}, Section 1.11). Furthermore, we can find a totally geodesic two-sphere $\tilde{\Sigma} \subset S^3$ such that $F(\tilde{\Sigma}) = \Sigma$. 

We next consider the induced map $i_\#: \pi_1(\Sigma) \to \pi_1(M)$. By Proposition \ref{incompressibility}, the map $i_\#$ is injective. We claim that $i_\#$ is surjective. To prove this, we consider a closed curve $\alpha: [0,1] \to M$. The path $\alpha$ induces an isometry $\psi: S^3 \to S^3$ satisfying $F \circ \psi = F$. Since $\tilde{\Sigma}$ and $\psi^{-1}(\tilde{\Sigma})$ are totally geodesic, we have $\tilde{\Sigma} \cap \psi^{-1}(\tilde{\Sigma}) \neq \emptyset$. Let us fix a point $\tilde{p} \in \tilde{\Sigma} \cap \psi^{-1}(\tilde{\Sigma})$. We can find a smooth path $\tilde{\gamma}: [0,1] \to \tilde{\Sigma}$ such that $\tilde{\gamma}(0) = \tilde{p}$ and $\tilde{\gamma}(1) = \psi(\tilde{p})$. We next define a smooth path $\gamma: [0,1] \to \Sigma$ by $\gamma(s) = F(\tilde{\gamma}(s))$. Clearly, $\gamma$ is a closed curve in $\Sigma$, i.e. $\gamma(0) = \gamma(1)$. Furthermore, $\gamma$ is homotopic to $\alpha$. Thus, we conclude that $[\alpha] = [\gamma] \in i_\#(\pi_1(\Sigma))$. This shows that the map $i_\#: \pi_1(\Sigma) \to \pi_1(M)$ is surjective.

\end{document}